\documentclass[11pt]{elsarticle}

\usepackage{amsmath, amssymb, color}

\newtheorem{thm}{Theorem}
\newtheorem{fact}[thm]{Fact}
\newtheorem{conjecture}[thm]{Conjecture}

\newtheorem{corollary}[thm]{Corollary}
\newtheorem{cor}[thm]{Corollary}
\newtheorem{lem}[thm]{Lemma}

\newtheorem{rem}[thm]{Remark}
\def\PP{{\mathcal P}}
\def\P{{\mathbb P}}
\def\E{{\mathbb E}}

\let\epsilon=\varepsilon

\newcommand{\eop}{\hfill$\Box$}
\newcommand{\textdef}{\textit}

\newcommand{\psn}{{\rm psn}}

\begin{document}
\begin{frontmatter}
\title{On the path separation number of graphs}

\author[Balogh]{J\'ozsef Balogh\fnref{GrantBal}}
\ead{jobal@math.uiuc.edu}

\author[Csaba]{B\'ela Csaba\fnref{GrantCsP,GrantCs}}
\ead{bcsaba@math.u-szeged.hu}

\author[Martin]{Ryan R. Martin\corref{cor}\fnref{ISU,GrantM}}
\ead{rymartin@iastate.edu}

\author[Pluhar]{Andr\'as Pluh\'ar\fnref{Szeged,GrantCsP}}
\ead{pluhar@inf.u-szeged.hu}

\address[Balogh]{Department of Mathematical Sciences, University of Illinois and Bolyai Institute, University of Szeged}
\address[Csaba]{Bolyai Institute, University of Szeged}
\address[Martin]{Department of Mathematics, Iowa State University, Ames, Iowa 50011}
\address[Pluhar]{Department of Computer Science, University of Szeged}

\fntext[GrantBal]{The first author's research is partially supported by NSF CAREER Grant DMS-0745185,
T\'AMOP-4.2.1/B-09/1/KONV-2010-0005, and Marie Curie FP7-PEOPLE-2012-IIF 327763.}
\fntext[GrantCsP]{The second and fourth authors were partially supported by
the European Union and the European Social Fund through project
FuturICT.hu (grant no.: TAMOP-4.2.2.C-11/1/KONV-2012-0013).}
\fntext[GrantCs]{The second author was also supported in part by the ERC-AdG. 321104.}
\fntext[GrantM]{The third author's research partially was supported by NSF grant DMS-0901008 and by NSA grant H98230-13-1-0226.}

\cortext[cor]{Corresponding Author. Ph: +1 515 294-1282. Fax: +1 515 294 5454.}

\begin{abstract}
A {\it path separator} of a graph $G$ is a set of paths $\mathcal{P}=\{P_1,\ldots,P_t\}$ such that for every pair of edges $e,f\in E(G)$, there exist paths $P_e,P_f\in\mathcal{P}$ such that $e\in E(P_e)$, $f\not\in E(P_e)$, $e\not\in E(P_f)$ and $f\in E(P_f)$. The {\it path separation number} of $G$, denoted ${\rm psn}(G)$, is the smallest number of paths in a path separator. We shall estimate the path separation number of several graph families -- including complete graphs, random graph, the hypercube -- and discuss general graphs as well.
\end{abstract}

\begin{keyword}
network reliability \sep test sets \sep path covering \sep path separation \sep path decomposition \sep trees \sep \MSC[2010]{05C35, 05C70}
\end{keyword}

\end{frontmatter}

\section{Introduction}
\label{sec:intro}
\thispagestyle{empty}

Separation of combinatorial objects is a well-studied question in mathematics and engineering, dating back to at least R\'enyi~\cite{R}.  In fact the notion goes by many terms: identification or, in engineering, {\em testing} are also used for the same idea~\cite{BS, FKKLN, FK2, HPWYC, HKL, TWH}.

Let $\mathcal{H}=(X, E)$ be a hypergraph with ground set $X$ and edge set $E$. We say that $L \subset E$ is a {\it weak separating system} if for all $x, y \in X,$ $x\neq y$ there exists an $A \in L$ such that either $x \in A$ or $y \in A,$ but $\{x, y\} \not \subset A$. Similarly, $L$ is a {\it strong (or complete) separating system} if for all $x, y \in X,$ $x\neq y$ there exist $A_x, A_y \in L$ such $x \in A_x$ and $y \in A_y$, but $x \not \in A_y$ and $y \not \in A_x$, as defined by Dickson~\cite{D}. Observe that any strong separating system is also a weak separating system. In several applications $X$ is just the vertices or edges of a certain graph $G$, while $E$ can be a set of closed neighborhoods, cycles, closed walks, paths, etc. of $G$, see e.g.~\cite{FGNPV,  FKKLN, HKL, TWH}.

In this paper we consider strong separation of the edges of graphs by paths. Since we deal with strong separation in this paper, we will just use the term ``separating system'' or ``separator'' when referring to a strong separating system. Let $G$ be a graph with at least two edges. A \textdef{path separator}
of $G$ is a set of paths $\mathcal{P}=\{P_1,\ldots,P_t\}$ such that for every pair of distinct edges $e,f\in E(G)$, there exist paths $P_e, P_f\in\mathcal{S}$ such that $e\in E(P_e)$ and $f\in E(P_f)$ but
$e\not\in E(P_f)$ and $f\not\in E(P_e)$. The \textdef{path separation number} of $G$, denoted $\psn(G)$, is the smallest number of paths in a path separator. If $G$ has exactly one edge then we say that $\psn(G)=1$ and if $G$ is empty then we say that $\psn(G)=0$.

R\'enyi~\cite{R} conjectured that $O(n)$ paths suffice for the weak separation in any graph with $n$ vertices. This problem is still unsolved, although Falgas-Ravry, Kittipassorn, Kor\'andi, Letzer and Narayanan~\cite{FGNPV} recently made some progress for proving it for trees and certain random graphs. We propose the stronger Conjecture~\ref{conj:bigo} below: $O(n)$ paths are sufficient even for strong separation.

In this paper we prove this conjecture for complete graphs (Theorem~\ref{thm:complete}), forests (Theorem~\ref{thm:forest}), higher dimensional cubes (Theorem~\ref{thm:hypercube}), and not too sparse random graphs (Theorem~\ref{thm:Gnp}). It is somehow surprising since generally strong separation may need many more paths than weak separation, as we remark following Theorem~\ref{thm:hypercube}.

Denote  $H_2(x)$ to be the {\it binary entropy function}, i.e.~$H_2(x)=-x\log_2 x - (1-x)\log_2(1-x)$, where $x\in (0,1)$. Denote $K_n$ to be the complete graph and $P_n$ to be the path  on $n$ vertices.
The parameters $\delta(G)$ and $\Delta(G)$ denote the minimum and maximum degree of $G$, respectively.

Fact~\ref{fact:edgenumber} follows from the fact that the edge set itself is a path separator if there are at least 2 edges and $\psn(G)=m$ if $m=1$ or $m=0$.

\begin{fact}\label{fact:edgenumber}
Let $G$ be a graph with $m$ edges. Then $\psn(G) \leq m$.
\end{fact}

Because of Fact~\ref{fact:sum} below, we will always assume that the graph $G$ that we are working with is connected.

\begin{fact}\label{fact:sum}
If $G$ is a graph that is the vertex-disjoint union of graphs $G_1$ and $G_2$ then $\psn(G)=\psn(G_1)+\psn(G_2)$.
\end{fact}

When $G$ is a forest we determine $\psn(G)$ in Theorem~\ref{thm:forest}, otherwise Theorem~\ref{thm:bounds} estimates it. Note that the proof of the lower bound  in Theorem~\ref{thm:bounds} does not use the structure of paths, only that a path has at most $n-1$ edges.

%
%

\begin{thm}\label{thm:bounds}
Let $G$ be a  graph on $n\geq 4$ vertices and $m\ge 2(n-1)$ edges, then
$$ \frac{m\ln m}{n\ln (en/2)} < \frac{\log_2 m}{H_2\left((n-1)/m\right)} \leq \psn(G) \leq 4n \lceil \log_2\lceil m/n\rceil\rceil +2n < 5n\log_2n.$$
\end{thm}

Theorem~\ref{thm:complete} establishes that the path separation number of the complete graph is at most $2n+4$ and  Theorem~\ref{thm:bounds} implies that it is at least $(1-o(1))n$. Of course, the bound $2n+4$ is not sharp even for $n=5$ or $6$, since by Fact~\ref{fact:edgenumber} we have that $\psn(K_n) \leq n(n-1)/2 <2n+4$ in these cases.

%


%

\begin{thm}\label{thm:complete}
For $n\ge 10$ we have $\psn(K_n)\le 4\lceil n/2\rceil +2\le 2n+4$.
\end{thm}

Theorem~\ref{thm:forest} gives an explicit formula for the path separation number of a forest $F$ depending only on the degree sequence and the number of connected components of $F$ that are, themselves, paths.  A
\textdef{path-component} of a graph is a connected component that is a path.

\begin{thm}\label{thm:forest}
Let $F$ be a forest with $v_1$ vertices of degree 1, $v_2$ vertices of degree 2 and $p$ path-components.  Then $\psn(F)=v_1+v_2-p$.
\end{thm}

\begin{corollary}\label{mintree}
The smallest path separation number for a tree $T$ on $n$ vertices is $\lceil n/2\rceil+1$. This is achieved with equality if and only if (a) $n$ is even and all the degrees of $T$ are either $1$ or $3$ or (b) $n$ is odd, $T$ has one vertex of degree either $2$ or $4$ and all other vertices have degree either $1$ or $3$.
\end{corollary}

\begin{corollary}\label{star}
If $G$ is a tree with $n$ vertices then $psn(G)=n-1$ if and only if $G$ is a subdivided star.
\end{corollary}

Theorem~\ref{thm:hypercube} considers the graph of the $d$-dimensional hypercube $Q_d$, whose path separation number shows different behavior from our previous results.

\begin{thm}\label{thm:hypercube}
For $d\geq 2$, let $Q_d$ denote the $d$-dimensional hypercube. Then $\frac{d^2}{4\ln d}\leq\psn(Q_d)\leq 3d^2+d-4$.
\end{thm}

Theorem~\ref{thm:hypercube} also demonstrates the difference between weak and strong separation: Honkala, Karpovsky and Litsyn proved in~\cite{HKL} that essentially $d+\log_2 d$ cycles are necessary and sufficient for a weak separation of the edges of the hypercube, which easily translates to a weak path separator having essentially at most $2(d+\log_2 d)$ paths, that is, much less than what is required in any strong separating system.


In Theorem~\ref{thm:Gnp} we address the Erd\H{o}s-R\'enyi random graph in which each pair of vertices is, independently, chosen to be an edge with probability $p$. We say that a sequence of random events occurs \textit{with high probability} if the probability of the events approaches 1 as $n\rightarrow\infty$.

\begin{thm}\label{thm:Gnp}
Let  $p=p(n)>1000\log n/n$ and $s=4\log n/\log(pn/\log n)$. Then $\psn(G(n,p)) = O(psn)$ with high probability.

In particular, for $\alpha>0$ and $p=p(n)>n^{\alpha-1}$ this gives $\psn(G(n,p))=\Theta(pn)$ with high probability and
for $p=p(n)>10\log n/n$ it yields that $\psn(G(n,p))=O(pn\log n)$, with high probability.
\end{thm}



Although Theorem~\ref{thm:bounds} establishes that $\psn(G)=O(n\log n)$ for any graph $G$ on $n$ vertices, it is not clear that the path separation number of a graph is that large. We can ensure that there are graphs $G$ on $n$ vertices with $\psn(G)\geq (2-o(1))n$, which is larger than the lower bound in Theorem~\ref{thm:complete}. We prove Theorem~\ref{thm:bipartite} as a remark in the proof of Theorem~\ref{thm:bounds} because the techniques are similar.
\begin{thm}\label{thm:bipartite}
Let $K_{a,b}$ denote the complete bipartite graph with $a$ vertices in one part and $b$ in the other. Fix $\epsilon\in (0,1/2)$. For $n$ sufficiently large, $\psn(K_{\epsilon n,(1-\epsilon)n})>2(1-2\epsilon)n$.
\end{thm}

In the rest of the paper we shall prove the above theorems. We close Section~\ref{sec:intro} with the following conjecture (a weaker version of it, for weakly separating systems, was formulated in~\cite{FKKLN}).

\begin{conjecture}\label{conj:bigo}
There exists a constant $C$ such that for every positive integer $n$ and for every graph $G$ on $n$ vertices $\psn(G)\le Cn.$
\end{conjecture}

Of course, since Fact~\ref{fact:edgenumber} gives an upper bound of $m$ on the path separation number, Conjecture~\ref{conj:bigo} is satisfied for any graph with $O(n)$ edges.  From Theorem~\ref{thm:bipartite}, we know that there are graphs that have path separation number arbitrarily close to $2n$.

\section{Proofs}
\label{sec:proofs}

\subsection{Proof of Theorem~\ref{thm:bounds}:}
We note that there is a trivial lower bound of $\lceil\log_2 m\rceil$ which follows from the fact that if $\mathcal{P}=\{P_1,\ldots,P_t\}$ is a path separator, then each edge has a different nonzero indicator vector $(X_1,\ldots,X_t)$ where $X_i=1$ if the edge is in $P_i$ and $0$ otherwise.  Since $H_2(x)\leq 1$, the lower bound in Theorem~\ref{thm:bounds} strictly improves on this trivial bound.~\\

\noindent\textbf{Proof of Theorem~\ref{thm:bounds}: Lower bound.} We use the entropy method. For facts about the entropy method, see Section 22 of Jukna~\cite{J}. The entropy of discrete random variable $Y$ is $H_2(Y):=\sum_i -\Pr(Y=y_i)\log_2\Pr(Y=y_i)$, where $\{y_i\}_i$ is the range of values of $Y$. The notation $H_2(p)$ denotes $-p\log_2 p-(1-p)\log_2(1-p)$ if $p$ is a real number in $(0,1)$ and $H_2(Y)$ denotes the entropy of random variable $Y$. The use will be clear from the context. Note also that if $X$ is a Bernoulli$(p)$ random variable, then $H_2(X)=H_2(p)$.

If $Y$ is a random variable that takes on $m$ values, each with equal probability, then $H_2(Y)=\log_2m$.  The subadditivity property of entropy says that if $Y$ can be expressed as the ordered tuple of random variables $(X_1,\ldots,X_t)$, then $H_2(Y)\leq\sum_{i=1}^tH_2(X_i)$.

Let $\pi_1,\ldots,\pi_t$
be the paths of a path separator of graph $G$. Let $X_i$ be the event that
a randomly-chosen edge is in path $\pi_i$.  Since the joint distribution
$(X_1,\ldots,X_t)$ takes on $m$ values each of which being equally likely, $H_2(X_1,\ldots,X_t)=\log_2m$.
Using the subadditivity property,
\begin{align}\label{eq:subadditivity}
  \log_2m &= H_2(X_1,\ldots,X_t)
  \leq \sum_{i=1}^tH_2(X_i)
  = \sum_{i=1}^tH_2\left(\frac{\ell(\pi_i)}{m}\right)
  \leq t H_2\left(\frac{n-1}{m}\right),
\end{align}
because every path has length at most $n-1$. In the last inequality we also used the fact that $H_2(p)$ is increasing for $p\in (0,1/2)$ together with the condition $m\ge 2(n-1).$

Writing
$x=({n-1})/{m}$ we have
\begin{align*}
  t &
\geq\frac{\log_2 m}{H_2\left(\frac{n-1}{m}\right)}
\geq \frac{\ln m}{-x\ln x-(1-x)\ln (1-x)}
  > \frac{\ln m}{-x\ln x+x} \\
  &= \frac{m\ln m}{n-1}\cdot\frac{1}{\ln\left(\frac{m}{n-1}\right)+1}
  \geq \frac{m\ln m}{n-1}\cdot\frac{1}{\ln\left(\frac{n}{2}\right)+1}
  = \frac{m\ln m}{(n-1)\ln (en/2)} .
\end{align*}

\begin{rem}[Proof of Theorem~\ref{thm:bipartite}]
   For the proof of this result, we use (\ref{eq:subadditivity}) to conclude that $t>\frac{\ln m}{-x\ln x+x}$ where $m=\epsilon(1-\epsilon)n^2$ and $x=(\epsilon n+1)/m$ because the length of any path in $K_{\epsilon n,(1-\epsilon)n}$ is at most $\epsilon n+1$.  For fixed $\epsilon$, the fraction $\frac{\ln m}{-x\ln x+x}$ goes to $2(1-\epsilon)n$ as $n\rightarrow\infty$. Thus for $n$ large enough, the lower bound $2(1-2\epsilon)n$ suffices.
\end{rem}

\noindent\textbf{Proof of Theorem~\ref{thm:bounds}: Upper bound.} We use a classical theorem of Lov\'asz~\cite{L}:
\begin{thm}[Lov\'asz~\cite{L}]\label{thm:Lovasz}
  The edges of a graph on $n$ vertices can be covered by at most
$\left\lfloor\frac{n}{2}\right\rfloor$ edge-disjoint paths and cycles.
\end{thm}

Since any cycle can be partitioned into two paths we obtain  the following corollary.
\begin{cor}\label{cor:Lovasz}
The edges of a graph on $n$ vertices can be covered by at most $n$ edge-disjoint paths.
\end{cor}

Returning to Theorem~\ref{thm:bounds}, apply Corollary~\ref{cor:Lovasz} to $G$ and partition $E(G)$ into at most $n$ paths, $P_1', \ldots, P_k'$, $k \leq n$. We shall cut each path that is longer than $m/n$ into paths of length $\lceil m/n\rceil$ and  possibly one shorter path.
The number of such paths is
$$ \sum_{i=1}^k\left\lceil\frac{|P_i'|}{\lceil m/n\rceil}\right\rceil \leq\sum_{i=1}^k\left\lceil\frac{n|P_i'|}{m}\right\rceil \leq\sum_{i=1}^k\left(\frac{n|P_i'|}{m}+1\right) \leq 2n $$
where $|P_i|$ is the length (number of edges) of $P_i$.
So we have that the new family $\mathcal{P}$ consists of at most $2n$ paths.

For each path $P\in \mathcal{P}$ label the edges of $P$  by binary vectors. Specifically, take an injective function $b_P: E(P) \rightarrow \{0,1\}^t$, where $t={\left\lceil \log_2 \lceil m/n\rceil\right\rceil}$. For each $i\in [t]$ define the graph $G_i$ whose vertex set is $V(G)$ and whose edge set consists of edges of each path $P\in \mathcal{P}$ whose $i$-th digit in its vector is $1$. We apply Corollary~\ref{cor:Lovasz} to each $G_i$, giving a path partition $\mathcal{P}_i$ of cardinality at most $2n$. Analogously, for each $i\in [t]$ define the graph $G'_i$ whose vertex set is $V(G)$ and whose edge set consists of edges of path $P\in \mathcal{P}$ whose $i$-th digit in its vector is $0$. Applying Corollary~\ref{cor:Lovasz} to each $G_i'$ gives a path partition $\mathcal{P}'_i$ of cardinality at most $2n$. Observe that $\cup_i \mathcal{P}_i \cup_i \mathcal{P}'_i \cup \mathcal{P}$ is of size at most $4n\left\lceil\log_2\lceil m/n\rceil\right\rceil+2n$.

Now we show that $\cup_i \mathcal{P}_i \cup_i \mathcal{P}'_i \cup \mathcal{P}$ is a path-separating system. Since $\mathcal{P}$ is a path-cover, it suffices to just consider edges $e,f$ that are in the same $P\in\mathcal{P}$. As such, the map $b_P$ ensures that the edges will differ in one coordinate, say $i$, and so without loss of generality $e\in E(G_i)-E(G_i')$ and $f\in E(G_i')-E(G_i)$. Thus, the path covers of $G_i$ and $G_i'$ suffice to separate $e$ and $f$.

The final inequality comes from using $m\leq\binom{n}{2}$. As such, $n\geq 4$ gives
\begin{align*}
   4n\left\lceil\log_2\lceil m/n\rceil\right\rceil+2n \leq 4n\left\lceil\log_2\lceil (n-1)/2\rceil\right\rceil+2n < 5n\log_2 n .
\end{align*}
This concludes the proof of Theorem~\ref{thm:bounds}.\eop

\subsection{Proof of Theorem~\ref{thm:complete}}

We start with a corollary of Theorem~\ref{thm:Lovasz}: The edges of $K_n$ can be covered by at most
$\lceil\frac{n}{2}\rceil$ edge-disjoint paths.

Note that Lov\'asz~\cite{L} proved this when $n$ is even. For $n$ odd, Theorem~\ref{thm:Lovasz} implies the existence of a partition of $K_n$ into $\lfloor n/2 \rfloor$ Hamiltonian cycles but one can check that it is always possible to choose one edge from each Hamilton cycle so that these edges could be covered with one additional path.

Denote such a collection of $\lceil n/2\rceil$ edge-disjoint paths by $\mathcal{P}_1=\{P_1,\ldots,
P_{\left\lceil n/2\right\rceil}\}$. Select three random permutations $\alpha$, $\beta$ and $\gamma$ uniformly and independently of each other from $S_n$, the set of $n$-element permutations. Let $\mathcal{P}_\alpha$, $\mathcal{P}_\beta$ and $\mathcal{P}_\gamma$ be the images of $\mathcal{P}_1$ under the permutations $\alpha, \beta$ and  $\gamma$, respectively. That is, if $P_i=\{i_1, \ldots, i_n\} \in \mathcal{P}_1$ then, say, $\alpha(P_1)=\{\alpha(i_1), \alpha(i_2), \dots, \alpha(i_n)\}$, assuming $V(K_n)=[n]$.

If a pair of edges $e,f$ are in different paths in $\mathcal{P}_1$ then they are separated by $\mathcal{P}_1$.
Otherwise the probability that they are {\em not} separated by
$\mathcal{P}_\alpha$ is at most $2/(n-2)$. The separations
by $\mathcal{P}_\alpha, \mathcal{P}_\beta$ and $\mathcal{P}_\gamma$ are
independent of each other, i.e. the probability that $e$ and
$f$ are not separated by the system of paths
$$ \mathcal{P}=\mathcal{P}_1 \cup \mathcal{P}_\alpha \cup \mathcal{P}_\beta \cup \mathcal{P}_\gamma $$
is at most $(2/(n-2))^3$. The number of pairs that are not separated by
$\mathcal{P}_1$ is at most $\binom{n-1}{2}\left\lceil\frac{n}{2}\right\rceil$ and so the
expected number of pairs not separated by
$\mathcal{P}$ is less than $(2/(n-2))^3\binom{n-1}{2}\left\lceil\frac{n}{2}\right\rceil < 3$ for $n\ge 10$.
We can finish the path separator with two additional paths that separate the remaining two pairs of edges.\eop

\subsection{Proof of Theorem~\ref{thm:forest}}
First we prove that  $\psn(P_k)=k-1$, where $P_k$ is a path with $k$ vertices $x_1,\ldots,x_k$ and edge set $\{x_1x_2,x_2x_3,\ldots,x_{k-1}x_k\}$. For the upper bound observe that the edge set itself is a path separator of size $k-1$. For the lower bound, one observes that $\{x_{k-1},x_k\}$ must be a member of any path separator, otherwise the edge $x_{k-1}x_k$ cannot be separated from $x_{k-2}x_{k-1}$. Any path separator of $\{x_1,\ldots,x_{k-1}\}$ has size at least $k-2$ and the lower bound follows. Let $v_i$ be number of vertices of degree $i$. By Fact~\ref{fact:sum} it is sufficient to prove only that $\psn(T)=v_1+v_2$ for any tree $T$, provided $T$ is different from a path.

Let us start with the upper bound. Let $T$ be a non-path tree on $n\geq 4$ vertices with $v_1$ leaves and $v_2$ vertices of degree $2$. Add an extra vertex (we may call it $\infty$) that is incident to all of the leaves of $T$. Call the new graph $T^{\infty}$. This graph is planar and we can form a path for every face of the embedding of $T^{\infty}$ into the plane that begins and ends at a leaf and contains every vertex (other than ``$\infty$'') of the face. Each leaf will be the endpoint of exactly two such paths. For each vertex $w$ of degree 2, take  one of these paths that passes through it and partition it into two subpaths (that is subgraphs that are themselves paths), both of which have $w$ as an endpoint. We continue this process until we have separated all edges incident to the same degree-$2$ vertex. In other words, the degree-$2$ vertices are `cut points' of these paths.

The size of this family of paths is $v_1+v_2$. This family clearly separates any pair of edges that are not on the same two faces (each edge is on exactly two faces). If two edges, $e$ and $f$, share the same pair of faces then there must be a path of degree 2 vertices (possibly without any additional edges, but containing at least one vertex) connecting them. But by the partition of the paths at vertices of degree 2, one of the original paths containing $e$ and $f$ was partitioned into at least two paths, one that contains $e$ but not $f$ and another that contains $f$ but not $e$.

For  the lower bound, we proceed by induction on the order of  $T$. For the base case, the smallest tree which is not a path is a star on $4$ vertices; easily a separating system of size $3$ exists. If a leaf edge is covered by only one path then this path necessarily has length 1. If there is a leaf that is covered by a path of length $1$ then both the leaf from the tree and the path from the system can be removed. Assume that there is a tree $T$ with a minimum-sized path separator $\{P_1,\ldots,P_t\}$ such that no leaf is in exactly one path. In this case we use a simple discharging argument. Let us give each $P_i$ a charge of $1$ and discharge $1/2$ to each of its endpoints in $T$.  Every degree-$2$ vertex $w$ must receive a charge of at least $1$, otherwise the incident edges of $w$ are not separated. Every leaf $x$ must receive a charge of at
least $1$ because there are at least $2$ paths that contain the edge incident to $x$. Thus the number of paths is at least  $v_1+v_2$.~\eop

\subsection{Proof of Corollary~\ref{mintree}}

Lemma~\ref{lem:leaves} with Theorem~\ref{thm:forest} implies that the smallest path separation number
for a tree on $n$ vertices is $\lceil n/2\rceil+1$.

\begin{lem}\label{lem:leaves}
Let $T$ be a tree on $n\geq 2$ vertices with $v_1$ vertices of degree
1 and $v_2$ vertices of degree 2. Then $v_1+v_2\geq\lceil n/2\rceil+1$.
This is achieved with equality if and only if (a) $n$ is even and all
the degrees of $T$ are either $1$ or $3$ or (b) $n$ is odd, then $T$ has
one vertex of degree either $2$ or $4$ and all other vertices have
degree either $1$ or $3$.
\end{lem}

\subsection{Proof of Lemma~\ref{lem:leaves}}
Let us use the notation that $v_i$ denotes the number of vertices of degree $i$ in $T$ in which case $\sum_i v_i=n$ and $\sum_i iv_i=2n-2$. This gives $3\sum_iv_i-\sum_iiv_i=n+2$. Rearranging,
\begin{align}
   2v_1+2v_2=n+2+v_2+\sum_{i\geq 3}(i-3)v_i\geq n+2 . \label{eq:v1v2}
\end{align}
Hence, $\psn(T)=v_1+v_2\geq n/2+1$.

Now we investigate the case where $v_1+v_2=\lceil n/2\rceil+1$. If $n$ is even, then (\ref{eq:v1v2}) gives that $v_2+\sum_{i\geq 3}(i-3)v_i=0$, implying that $v_i=0$ for $i=2$ and $i\geq 4$. In addition, it gives $v_1=n/2+1$.

If $n$ is odd, then (\ref{eq:v1v2}) gives that $v_2+\sum_{i\geq 3}(i-3)v_i=1$, implying either that $v_2=1$, $v_1=(n+1)/2$, and $v_i=0$ for $i\geq 4$ or that $v_4=1$, $v_1=(n+3)/2$, and $v_i=0$ for $i=2$ and $i\geq 5$. \eop

\subsection{Proof of Theorem~\ref{thm:hypercube}}
For the lower bound we use Theorem~\ref{thm:bounds}, noting that $Q_d$ has $2^d$ vertices and $d2^{d-1}$ edges:
\begin{equation}\label{eq:QdLB}
   \psn(Q_d)\ge \frac{\log_2(d 2^{d-1})}{H_2((2^d-1)/d2^{d-1})} .
\end{equation}

If $d=2$, then it is easy to see that $\psn(Q_2)=4>\frac{2^2}{4\ln 2}\approx 1.44$.

If $d=3$, then (\ref{eq:QdLB}) gives $\psn(Q_3)\geq \frac{\log_2(12)}{H_2(7/12)}>3.6>\frac{3^2}{4\ln 3}\approx 2.05$.

If $d\geq 4$, then $H_2((2^d-1)/d2^{d-1})\leq H_2(2/d)\leq \frac{2/d}{\ln 2}\left(1-\ln(2/d)\right)$. Therefore
\begin{align}
\psn(Q_d) &> \frac{\log_2(d 2^{d-1})}{H_2(2/d)} > \frac{\ln(d2^{d-1})}{(2/d)(1-\ln (2/d))} = \frac{d^2}{2\ln d}\left(\frac{\frac{d-1}{d}\ln 2+\frac{\ln d}{d}}{1+\frac{1-\ln 2}{\ln d}}\right) \label{eq:bigparen} \\
&> \frac{d^2}{2\ln d}\left(\frac{\ln 2}{1+\frac{1-\ln 2}{\ln 4}}\right) > \frac{d^2}{4\ln d} . \nonumber
\end{align}
The second inequality follows from $H_2(x)\geq x(1-\ln x)$, which holds for $x\leq 0.65$. The term in parentheses in (\ref{eq:bigparen}) is in fact larger than the lower bound of 1/2 that eventually results. However, it is easy to verify this bound and we have no interest in optimizing the constant.~\\

Consider the upper bound. We will show by induction on $d$ that $f(d)=3d^2+d-4$ paths suffice. As we have established above, $\psn(Q_2)=4$. It is helpful to view the vertices of $Q_d$ as $\{0,1\}$-vectors of dimension $d$.  Let $Q^i$ denote the $(d-1)$-dimensional subcube whose vertices have $i$ in the $d^{\rm \: th}$ coordinate for $i=0,1$.  Edges between vertices in $Q^i$ are called \textit{$i$-interior edges} for $i=0,1$. Edges with one endvertex in $Q^0$
and the other in $Q^1$ are called \textit{crossing edges}. Consider an edge $e^0$ in $Q^0$ and an edge $e^1$ in $Q^1$. We call such edges \textit{mirror images} if the endvertices of $e^0$ can be made into the endvertices of $e^1$ simply by changing the $d^{\rm \: th}$ coordinate from $0$ to $1$. For a path in $Q^0$ the \textit{mirror image path} is defined in an analogous way. We construct three different types of paths.

\textbf{Type 1:} By the inductive hypothesis, $Q^0$ has a path separation set of size $f(d-1)$. Construct it and its mirror image in $Q_1$.  Then for each pair of mirrored paths, connect their final endpoints via a crossing edge.  There are $f(d-1)$ paths of Type 1.

With this set of paths, the following pairs of edges $(e',e'')$ are separated: (1) if both are $0$-interior edges or both are $1$-interior edges or (2) if $e'$ is $0$-interior and $e''$ is $1$-interior but they are not mirror images.

So there are only three types of pairs of edges $(e',e'')$ that are not separated: (3) if $e'$ and $e''$ are mirror images or (4) if $e'$ is a crossing edge and $e''$ is an interior edge or (5) if both $e'$ and $e''$ are crossing edges.

\textbf{Type 2:} Construct an arbitrary system of paths that covers the edges of $Q^0$. It is easy to prove by induction that there is such a system consisting of $d$ paths. Construct the set of mirror image paths in $Q^1$. There are $2d$ paths of Type 2.

The set of Type 2 paths  separates pairs of mirror image edges (thus satisfying (3) above).  Furthermore, for every crossing edge $e'$ and interior edge $e''$ there is a path in this second group containing $e''$ but not $e'$. So they will be used to aid in separation of the pairs in (4).


\textbf{Type 3:} Construct a system of paths via a \textit{separating family} of the $2^{d-1}$ crossing edges of size $2d-2$.\footnote{A (strong) separating family, $F_0\cup F_1$, of a set $\Sigma$ of size $n$ can be found as follows: Assign a binary codes of length $\lceil\log_2 n\rceil$ to each member of $\Sigma$. Place it into the $j^{\rm th}$ member of $F_0$ if and only if the $j^{\rm th}$ bit of its code is $0$. Place it into the $j^{\rm th}$ member of $F_1$ if and only if the $j^{\rm th}$ bit of its code is $1$. Each of $F_0$ and $F_1$ is a weak separating family and their union is a strong separating family. The size of $F_0\cup F_1$ is $2\lceil\log_2 n\rceil$ and the sizes of the sets are at most $n/2$.} That is, sets $S_1,\ldots,S_{2d-2}$ of crossing edges so that for each pair of crossing edges there is an $S_i$ that contains the first but not the second and an $S_j$ that contains the second but not the first.

For each $S_j$ we will construct two paths utilizing a Hamilton path $v_1^0,\dots,v_{N}^0$ in the $(d-1)$-dimensional hypercube $Q^0$ and its mirror image in $Q^1$.  Here $N=2^{d-1}$. (Hypercubes of dimension at least 2 are well-known to be Hamiltonian. According to~\cite{Fink}, this fact goes back to 1872~\cite{Gros}.)

Let $S_j=\{v_{i_1}^0v_{i_1}^1,v_{i_2}^0v_{i_2}^1,\ldots,v_{i_N}^0v_{i_N}^1\}$ with $i_1<i_2<\cdots<i_N$. The first of the two paths related to $S_j$ will begin by traversing the Hamilton path in $Q^0$ crossing at the first opportunity and then traversing the mirror image Hamilton path in $Q^1$ crossing at the next opportunity and continuing in $Q^0$.  We continue this until we finish in $Q^0$:
$$ v_1^0,\ldots,v_{i_1}^0,v_{i_1}^1,v_{i_1+1}^1,\ldots,v_{i_2-1}^1,v_{i_2}^1,v_{i_2}^0,v_{i_2+1}^0,\ldots,v_{i_N-1}^1,v_{i_N}^1,v_{i_N}^0,v_{i_N+1}^0,\ldots,v_{2^{d-1}}^0 . $$
The second path is a mirror image:
$$ v_1^1,\ldots,v_{i_1}^1,v_{i_1}^0,v_{i_1+1}^0,\ldots,v_{i_2-1}^0,v_{i_2}^0,v_{i_2}^1,v_{i_2+1}^1,\ldots,v_{i_N-1}^0,v_{i_N}^0,v_{i_N}^1,v_{i_N+1}^1,\ldots,v_{2^{d-1}}^1 . $$
This group of paths is of size $2(2d-2)$ and separates pairs of crossing edges; that is, those pairs in (5).

To complete the separation of the paths in (4), we need that for every crossing edge $e'$ and interior edge $e''$ there is a path containing $e'$ but not $e''$. If we call $S_0$ the set of {\em all} crossing edges, we can construct two additional paths. Add two new (Hamiltonian) paths to our system, first starting the alternating path from an $x^0 \in Q^0$, then from $x^1 \in Q^1$. These paths contain every crossing edge $e'$ but every interior edge $e''$ is left out from at least one of those. There are $2(2d-2)+2=4d-2$ paths of Type 3.

The total number of paths in the three groups is $f(d-1)+2d+4d-2=f(d-1)+6d-2$. Setting $f(d)=3d^2+d-4$ satisfies $f(2)=4$ and $f(d)=f(d-1)+6d-2$. \eop~\\

\begin{rem} We believe that the lower bound is correct in the sense that $\psn(Q_d)=O(d^2/\log d)$. We think that a proof is likely in the same vein as the proof of the upper bound of $\psn(K_n)$. Note that $E(Q_d)$ can be covered by $d$ paths as we have described in the above proof.  Fix such a path system $\PP={P_1,\ldots,P_d}$. Now choose, randomly and independently, $100d/\log d$ automorphisms of $Q_d$ and apply it to $\PP$. This will give a path system $\PP^*$ of size $100d^2/\log d$. The system $\PP$ does not separate at most $d2^{2d}$ pairs of edges, the ones which are in the same path $P_i$ for some $i$. Unfortunately we do not have a good estimate on the probability that a pair of edges $(e',e'')$ are not separated in $\PP^*$, unless we know that no $P_i$ contains more than $O(d/\log d)$ edges that are crossing with respect to a given partition.
\end{rem}

\subsection{Proof of Theorem~\ref{thm:Gnp}}

The idea of the proof is to partition $G(n,p)$ into several random graphs such that every pair of edges should be separated by them and every pair of edges should be in several of the random graphs. Then by Vizing's theorem, we partition the edges of each of the random graphs into matchings and using some other random graphs we connect them into paths.

Let $G$ be a graph chosen according to the distribution $G(n,p)$. First we partition $E(G)$ into four random graphs. Let $f: E(G) \to \{1,2,3,4\}$ be a function so that each $f(e)$ is chosen uniformly from $\{1,2,3,4\}$ independently. Let $E_i:\{e: f(e)=i\}$ for $1\le i\le 4$.

Next we form six random subgraphs with the same vertex set $V(G)$ and with edge sets as follows: $E(G_1^1)=E_1\cup E_2$,  $E(G_1^2)=E_3\cup E_4$,  $E(G_2^1)=E_1\cup E_3$,  $E(G_2^2)=E_2\cup E_4$,  $E(G_3^1)=E_1\cup E_4$,  $E(G_3^2)=E_2\cup E_3$. These six random subgraphs have the property that for any pair of edges $e,f\in E(G)$ there is an $i,j$ that $e,f\in E(G_i^j)$.

Now fix a pair of indices $(i,j)$. Without loss of generality, assume $i=j=1$ and consider $G_1^1$. Note that $G_1^1$ is itself a random graph distributed according to $G(n,p/2)$.
We will further partition the edge-set of $G_1^1$ into random subgraphs.

Fix $r=\lfloor 3pn/(64\log n)\rfloor
>40$
and let $g: E(G^1_1) \to \{1,\ldots,r\}$ be a function so that each $g(e)$ is chosen uniformly from
$\{1,\ldots,r\}$ independently. Repeat this process $s=\lfloor 6\log n/\log(pn/\log n)\rfloor$ times. Because $p>1000\log n/n$,
we have $sr\leq\frac{6\log n}{\log(pn/\log n)}\cdot\frac{3pn}{64\log n}=O(n)$ subgraphs $H_1,\ldots,H_{sr}$, each of which is a copy
of $G(n,p/(2r))$.

The set of graphs $\{H_{\alpha} : \alpha=1,\ldots,sr\}$ will separate every pair of edges with high probability. To see this, the union bound gives
\begin{align*}
   {\P}\left(\exists e,f\in E(G_1^1) : \mbox{no $H_{\alpha}$ separates $e,f$}\right) &\leq \sum_{e,f\in E(G_1^1)}{\P}\left(\mbox{no $H_{\alpha}$ separates $e,f$}\right) \\
   &\leq \binom{\binom{n}{2}}{2}\left(\frac{1}{r}\right)^s \\
   &\leq \exp\{4\log n-s\log r\}=O\left(n^{-2}\right) .
\end{align*}

Furthermore, with high probability, all of the graphs $H_{\alpha}$ have maximum degree less than $2n\frac{p}{2r}
\leq 25\log n$, noting that the average degree is $n\frac{p}{2r}$. By a Chernoff bound (Page 12 of Bollob\'as~\cite{BB}) and the union bound
\begin{align*}
   {\P}(\exists\alpha : \Delta(H_{\alpha})\geq pn/r) &\leq \sum_\alpha \sum_{v\in V(H_\alpha)} {\P}(\deg(v)\geq pn/r)
   \leq srn\exp\left\{-\frac{t^2}{2(\mu+t/3)}\right\} ,
\end{align*}
where $\mu={\E}[\deg(v)]=n\frac{p}{2r}$ and $t=\frac{pn}{r}-\mu=\frac{pn}{2r}$. Hence,
\begin{align*}
   {\P}(\exists\alpha : \Delta(H_{\alpha})\geq pn/r) &\leq srn\exp\left\{-\frac{3pn}{16r}\right\}= O(n^2\exp\{-4\log n\})=O(n^{-2}) .
\end{align*}

The idea for the rest of the proof is that by Vizing's theorem with high probability, $H_{\alpha}$ can be partitioned into at
most $25\log n\geq \Delta(H_{\alpha})+1$ matchings and by using edges of another subgraph of $G^2_1$, we connect them into paths.
The total number of paths used will be at most
\begin{align}\label{eq:spn}
   (25\log n)\cdot s\cdot pn/\log n=25spn . 
\end{align}

We will need an additional notion. Let $D(n,p)$ be the oriented directed graph on $n$ vertices;
i.e. for every~\textit{ordered pair} $(x, y)$, we draw an edge with probability $p$ independently of each other.
McDiarmid showed in~\cite{M} that the probability of the existence of  a Hamiltonian cycle in $D(n,p)$ is not less
than the same in $G(n,p)$, that is
$$ {\P}(G(n,p){\rm \; is \; Hamiltonian})\leq{\P}(D(n,p){\rm \; is \; Hamiltonian}).$$

In order to lower bound the probability that a random graph is Hamiltonian, we need to investigate more carefully the results that establish the threshold for Hamiltonicity. Although the paper of Koml\'os and Szemer\'edi~\cite{KSz} provides a very precise threshold, the classical paper of P\'osa~\cite{P} suffices for our needs and provides a simpler argument. (The point at which the evolution of the random graph achieves Hamiltonicity was established independently by Ajtai, Koml\'os and Szemer\'edi~\cite{AKSz} and by Bollob\'as~\cite{B1984}. See also Chapter 8 of Bollob\'as~\cite{BB}.) 

P\'osa's argument can be slightly modified to obtain that if $p\geq (30+9\alpha)\ln n/n$, then ${\P}(G(n,p){\rm \; is \; Hamiltonian})\geq 1-n^{-\alpha}$. In particular,
\begin{equation}\label{eq:Ham}
   p\geq 100\ln n/n \qquad \Longrightarrow \qquad {\P}(G(n,p){\rm \; is \; Hamiltonian})\geq 1-O\left(n^{-7}\right) ,
\end{equation}

So apply Vizing's theorem and decompose, say $H_1$, into $\Delta_0+1:=\Delta(H_1)+1$ matchings, where of course $\Delta(H_1)\leq25\log n$ with high probability. Label these matchings as $M_1,\dots,M_{\Delta_0+1}$. Then for each $i\in\{1,\ldots,\Delta_0+1\}$ we form a separating matching system $M_i^1, \dots, M_i^{t}$ for each $M_i$, where $t=2\lceil \log_2(\Delta_0+1)\rceil$. I.e., for every pair of edges in $M_i$ there is a pair $\{M_i^{j_1},M_i^{j_2}\}$ such that one edge is in $M_i^{j_1}-M_i^{j_2}$ and the other is in $M_i^{j_2}-M_i^{j_1}$. A separating matching system of this size is possible, as we saw in the proof of Theorem~\ref{thm:hypercube}.

Now for each $i\in [\Delta_0+1]$ and $j\in [t]$ we find a path containing the edges of $M_i^j\subset E(G^1_1)$ connected by the edges of $G^2_1$. For this we define an auxiliary oriented directed graph $D$. The vertex set of $D$ will be the union of $M_i^j$ and $V(G)-V(M_i^j)$. The edge set will be defined as follows: First assign an arbitrary orientation to the edges of $M_i^j$. For (oriented) edges $(xy),(uv)\in M_i^j$ we have the edge $(xy)(uv)\in E(D)$ if $yu\in E(G^2_1)$ and the edge $(uv)(xy)\in E(D)$ if $vx\in E(G^2_1)$. For an (oriented) edge $xy\in M_i^j$ and $u\in V(G)-V(M_i^j)$ we have in $D$ the oriented edge $(xy)u$ if $yu\in E(G^2_1)$ and we have the oriented edge $u(xy)$ if $xu\in E(G^2_1)$. If $w,z\in V(G)-V(M_i^j)$, then both oriented edges $(wz)$ and $(zw)$ are in $E(D)$ if $wz\in E(G^2_1)$.

The key property of $D$ is that if it contains a Hamilton path, then  $M_i^j\cup E(G^2_1)$ will contain the desired path. If $m$ denotes the number of vertices in $D$, then $n/2\leq m\leq n$. Each directed edge in $D$, however, is present with probability $p/2$. To see this, observe that given any orientation of $M_i^j$, edges of the form $(xy)(uv)$, the form $(xy)u$, the form $y(uv)$ and the form $yu$ depend on the presence of $yu\in E(G^2_1)$. The probability that $yu$ is present in $E(G^2_1)$ is $p/2$. So, the edge probability in $D$ is the same as in $D(m,p/2)$, which is at least that of $G(m,p/2)$, therefore it contains a Hamilton path with probability at least $1-n^{-7}$ by (\ref{eq:Ham}).

As we see in (\ref{eq:spn}), since the total number of such paths is at most $25spn=O(n\log n)$, the union bound gives that all of these paths can be constructed with probability at least $1-O(\log n/n^6)$, which concludes the proof.\eop

\bigskip

{\bf Added in proof.} Around the time that we were finishing writing up our results, a similar paper appeared in the arXiv by Falgas-Ravry,  Kittipassorn,  Kor\'andi,  Letzter,  and  Narayanan~\cite{FKKLN}. Their work is independent from ours, and they consider a different separation: for each pair of edges they are interested to find a separating set which contains exactly one of them. In many of the cases (like trees), this leads to a different behavior. In some other cases, similar proof techniques as in our paper might be applied.~\\

{\bf Acknowledgements.} We appreciate the hard work of referees who uncovered some errors in the original manuscript.

\bibliographystyle{plain}{}
\bibliography{PathSeparation}

\end{document}